\documentstyle[12pt, amstex]{article}

\begin{document}
\newcommand{\bx}{\hfill\rule{.25cm}{.25cm}\medbreak}
\newtheorem{prop}{Proposition}[section]
\newtheorem{defn}{Definition}[section]
\newtheorem{lem}{Lemma}
\newtheorem{thm}{Theorem}[section]
\newtheorem{asum}{Assumption}
\newtheorem{cor}{Corollary}

\newcommand{\bR}{\overline{R}}
\newcommand{\eg}{\frak g}
\newcommand{\stu}{{\frak stu}_n(R, -, \boldsymbol{\gamma})}
\newcommand{\eij}{e_{ij}(a)}
\newcommand{\ejk}{e_{jk}(b)}
\newcommand{\Tij}{T_{ij}(a, b)}
\newcommand{\eik}{e_{ik}(ab)}
\newcommand{\ekl}{e_{kl}(b)}
\newcommand{\st}  {{\frak st}_n(R)}
\newcommand{\sth}  {\widehat{{\frak st}}_n(R)}
\newcommand{\stf} {{\frak st}_4(R)}
\newcommand{\stfh} {\widehat{{\frak st}}_4(R)}
\newcommand{\stft} {\widetilde{{\frak st}}_4(R)}
\newcommand{\stfs} {{\frak st}_4(R)^\sharp}
\newcommand{\stt} {{\frak st}_3(R)}
\newcommand{\stth} {\widehat{{\frak st}}_3(R)}
\newcommand{\sttt} {\widetilde{{\frak st}}_3(R)}
\newcommand{\stts} {{\frak st}_3(R)^\sharp}
\newcommand{\lij}{\gamma_i\gamma_j^{-1}}
\newcommand{\lji}{\gamma_j\gamma_i^{-1}}
\newcommand{\ljk}{\gamma_j\gamma_k^{-1}}
\newcommand{\gd}{\dot{\frak{g}}}

\centerline{\large \bf Universal coverings of Steinberg Lie
algebras of small characteristic}

\vspace{1em} \centerline{\bf Yun Gao and Shikui
Shang\footnotetext{2000 Mathematics Subject Classification:
17B55, 17B60.\\
Research of the first author was partially supported by NSERC of
Canada and Chinese Academy of Science.}}

\vspace{1.5em}
\begin{abstract}
It is well-known that the second homology group $H_2(\st)$ of the
Steinberg Lie algebra $\st$ is trivial when $n\geq 5$. In this
paper, we will work out $H_2(\st)$ explicitly for $n=3, 4$ which
are not necessarily trivial. Consequently, we obtained
$H_2(sl_n(R))$ for $n=3, 4$.
\end{abstract}

\vspace{0.5em}

\noindent{\bf Introduction}

 Steinberg Lie algebras $\st$ and/or
their universal coverings have been studied by Bloch [Bl],
Kassel-Loday [KL], Kassel [Ka], Faulkner [F], Allison-Faulkner
[AF], Berman-Moody [BM], [G1, 2] and [AG], and among others. They
are Lie algebras graded by finite root systems of type $A_l$ with
$l\geq 2$. In most situations, the Steinberg Lie algebra $\st$ is
the universal covering of the Lie algebra $sl_n(R)$ whose kernel
is isomorphic to the first cyclic homology group $HC_1(R)$ of the
associative algebra $R$ and the second Lie algebra homology group
$H_2(\st)=0$. It was shown in [Bl] and [KL] that $H_2(\st)=0$ for
$n\geq 5$.  [KL] mentioned without proof that $H_2(\st)=0$ for
$n=3, 4$ if $\frac{1}{2}\in$ lies in the base ring $K$. This was
proved (see [G1] 2.63) for $n=3$ if $\frac{1}{6}\in K$ and for
$n=4$ if $\frac{1}{2}\in K$.

In this paper, we shall work out $H_2(\st)$ explicitly for $n=3,
4$ without any assumption on (characteristic of) $K$. It is
equivalent to work on the Steinberg Lie algebras $\st$ of small
characteristic for small $n$. This completes the determination of
the universal coverings of the Lie algebras $\st$ and $sl_n(R)$ as
well.

More precisely, let $K$ be a unital commutative ring  and $R$ be a
unital associative $K$-algebra. Assume that $R$ has a $K$-basis
containing the identity element (so $R$ is a free $K$-module). The
Lie algebra $sl_n(R)$ is the subalgebra of  $gl_n(R)$ (the $n$ by
$n$ matrix Lie algebra over $K$ with coefficients in $R$),
generated by $\eij$ for $ 1\leq i\neq j \leq n$, $a\in R$, where
$e_{ij}$ is the standard matrix unit. The elements $\eij$ satisfy
certain canonical relations. The Steinberg Lie algebra $\st$ is
defined by generators corresponding to $\eij$ and those same
canonical relations.

Recall that the radical $r(n)$ of a positive integer
$n=p_1^{j_1}\cdots p_k^{j_k}$ is defined to be $p_1\cdots p_k$,
where $p_1, \dots, p_k$ are distinct prime numbers and $j_1,
\dots, j_k$ are positive integers. For example, $r(3)=3$ and
$r(4)=2$. For any positive integer $m$, set
$${\mathcal I}_m = mR + R[R, R] \ \ \  \text{ and } \ \ \ R_m =
R/{\mathcal I}_m.$$

Our main result of this paper is the following.

\newtheorem{mthm}{Theorem}
\renewcommand{\themthm}{}
\begin{mthm} let $K$ be a unital commutative ring  and $R$ be a
unital associative $K$-algebra. Assume that $R$ has a $K$-basis
containing the identity element. Then
$$H_2(\st)= {R_{r(n)}}^6,$$
the direct sum of six copies of $R_{r(n)}$, for $n=3, 4$. In
particular,  if $R$ is commutative and the characteristic of $K$
is $2$ for $n=4$ and $3$ for $n=3$, then
$$H_2(\st)= R^6$$
\end{mthm}

 The organization of this paper is as follows. In Section 1,
we review some basic facts on Steinberg Lie algebras $\st$.
Section 2 will treat the $n=4$ case. The novelty here is the use
of cosets of the Klein four subgroup  in the symmetric group $S_4$
which acts on $\stf$ naturally. Section 3 will handle the $n=3$
case. Finally in Section 4 we make a few concluding remarks.

\vspace{4mm}

\noindent{\bf \S 1 Basics on $\st$}

Let $K$ be a unital commutative ring and $R$ be a unital
associative $K$-algebra. We always assume that $R$ has a $K$-basis
$\{r_\lambda\}_{\lambda\in\Lambda}$($\Lambda$ is an index set),
which contains the identity element $1$ of $R$, i.e.
$1\in\{r_\lambda\}_{\lambda\in\Lambda}$.

The $K$-Lie algebra of $n\times n$ matrices with coefficients in
$R$ is denoted by $gl_n(R)$. For $n\geq 2$, the elementary Lie
algebra $sl_n(R)$ (or $e_n(R)$) is the subalgebra of $gl_n(R)$
generated by the elements $e_{ij}(a)$, $a\in R$, $1\leq i\neq
j\leq n$, where $e_{ij}$ are standard matrix units. Note that
$sl_n(R)$ can be equivalently defined as $sl_n(R) =[gl_n(R),
gl_n(R)]$, the derived subalgebra of $gl_n(R)$, or $sl_n(R)=\{X\in
gl_n(R)| \text{ tr}(X) \in [R, R]\}$.

Clearly, for any $a, b\in R$,
\begin{equation} [\eij , \ejk ]
= \eik  \notag \tag{1.1}\end{equation} if $i, j, k$ are distinct
and
 \begin{equation} [\eij , \ekl ] = 0 \tag{1.2}\end{equation}
if $j\neq k, i \neq l$.

For $n\geq 3$, the Steinberg Lie algebra $\st$ is defined to be
the Lie algebra over $K$ generated by the symbols $X_{ij}(a)$,
$a\in R$, $1\leq i\neq j\leq n$, subject to the relations

\noindent ( see [Ka], [F] or [BM]):
\begin{align} &a\mapsto X_{ij}(a) \text{ is a $K$-linear map,}\tag{1.3}\\
&[X_{ij}(a), X_{jk}(b)] = X_{ik}(ab), \text{ for distinct } i, j, k, \tag{1.4}\\
&[X_{ij}(a), X_{kl}(b)] = 0, \text{ for } j\neq k, i\neq l,
\tag{1.5}
\end{align}
where $ a, b\in R,$ $ 1\leq i, j, k, l \leq n$.

Both Lie algebras $sl_n(R)$ and $\st$ are perfect(a Lie algebra
$\eg$ over $K$ is called perfect if $[\eg , \eg ] = \eg $). The
Lie algebra epimorphism:
\begin{equation} \phi: \st \to sl_n(R), \tag{1.6}\end{equation}
such that $\phi(X_{ij}(a)) = \eij$, is a covering (or central
extension) and the kernel of $\phi$ is isomorphic to $HC_1(R)$,
which is the first cyclic homology group of $R$([KL] or [L]). So
the universal covering of $sl_n(R)$ is also the universal covering
of $\st$ denoted by $\widehat{\frak{st}}_n(R)$. Our purpose is to
calculate $\widehat{\frak{st}}_n(R)$ for any ring $K$ and $n\geq
3$.

The following proposition can be similarly proved as in [AF] for
the unitary case.
\newtheorem{fprop}{Lemma 1.7}
\renewcommand{\thefprop}{}
\begin{fprop} Let ${\frak T}: = \sum_{1\leq i < j\leq n}[ X_{ij}(R), X_{ji}(R)]$.
Then ${\frak T}$ is a subalgebra of $\st$ containing the center
$\frak Z$ of $\st$ with $[{\frak T}, X_{ij}(R)] \subseteq
X_{ij}(R)$. Moreover,
\begin{equation}\st = {\frak T}\oplus_{1\leq i \neq j\leq n}X_{ij}(R).\tag{1.8}
\end{equation}
\end{fprop}

As for the decomposition of $\st$, we take
$\{r_\lambda\}_{\lambda\in\Lambda}$, the fixed $K$-basis of $R$,
then $\{X_{ij}(r)\}$
$({r\in\{r_\lambda\}_{\lambda\in\Lambda},1\leq i\neq j\leq n})$
can be extended to a $K$-basis $\Gamma$ of $\st$.

In fact, the subalgebra $\frak T$ has a more refined structure.
Setting
\begin{equation} \Tij = [ X_{ij}(a), X_{ji}(b)],
\tag{1.9}\end{equation}
\begin{equation} t(a, b) = T_{1j}(a, b) - T_{1j}(1, ba),\tag{1.10}
\end{equation} for $a, b\in R, 1\leq i\neq j\leq n$. Then t(a,b)
does not depend on the choices of $j$(see [KL]). Note that $\Tij$
is $K$-bilinear, and so is $t(a,b)$.

One  can  easily prove the following lemma (see [KL] or [G1]).
\newtheorem{blem}{Lemma 1.11}
\renewcommand{\theblem}{}
\begin{blem}
Every element $x\in \frak{T}$ can be written as
$$ x = \sum_{i}t(a_i, b_i) + \sum_{2\leq j\leq n} T_{1j}(1, c_j),$$
where $a_i, b_i, c_j\in R$.
\end{blem}

The following result is well-known(see [Bl] and [KL]).
\newtheorem{blthm}{Theorem 1.12}
\renewcommand{\theblthm}{}
\begin{blthm} {\normalshape} If $n\geq 5$,  then $\phi: {\frak st}_n(R) \to sl_n(R)$ gives the universal covering
of $sl_n(R)$ and so the second homology group of Lie algebra $\st$
 is $H_2(\st)= 0$.
\end{blthm}

For later use, we collect some formulas as follows. They can be
proved by using the Jacobi identity, see [KL], [G1] or [AG].
\begin{align} & T_{ij}(a, b)= -T_{ji}(b, a)\notag \\
&[T_{ij}(a, b), X_{kl}(c)]=0 \text{ for distinct }i, j, k, l\notag\\
&[T_{ij}(a, b), X_{ik}(c)]= X_{ik}(abc), \ \ [T_{ij}(a, b), X_{ki}(c)]= -X_{ki}(cab) \notag\\
&[T_{ij}(a, b), X_{ij}(c)]= X_{ij}(abc + cba) \notag \\
&[t(a, b), X_{1i}(c)]= X_{1i}((ab-ba)c), \ \ [t(a, b), X_{i1}(c)] = -X_{i1}(c(ab-ba)) \notag\\
 &[t(a, b), X_{jk}(c)]=0
\text{ for }j, k\geq 2\tag{1.13}
\end{align}

\vspace{4mm}

\noindent{\bf \S 2 Coverings of $\stf$}

In this section, we  compute the universal covering $\stfh$ of
$\stf$. In fact, if $K$ is a field and char $k\neq2$, $\stf$ is
central closed and the universal covering is itself (see
[G1,Corollary 2.63]).

Now we don't put any assumption on the characteristic of $K$.

For any positive integer $m$, let ${\mathcal I}_m$ be the ideal of
$R$ generated by the elements: $ma$ and $ab-ba$, for $a,b \in R$.
Immediately, we have
\newtheorem{flem}{Lemma 2.1}
\renewcommand{\theflem}{}
\begin{flem}
${\mathcal I}_m=mR+R[R,R] \text{ and } [R, R]R=[R, R]R$
\end{flem}
{\bf Proof:} Since $ma$ and $ab-ba$ generate $\mathcal I$, then
$${\mathcal I}_m=mR+R[R,R]+[R,R]R+R[R,R]R.$$
Since
$$[R,R]R\subseteq R[R,R]+[[R,R],R]\subseteq R[R,R]+[R,R]\subseteq R[R,R]$$
and similarly, $R[R, R]\subseteq [R, R]R$, $[R, R]R=[R, R]R$ and
so the lemma is proved.  $\Box$

Let $$R_m:=R/{\mathcal I}_m$$ be the quotient algebra over $K$
which is commutative.  Write $\bar a=a+{\mathcal I}_m$ for $a\in
R$. Note that if $m=2$ then  $\overline{a}=-\overline{a}$ in
$R_m$.

\newtheorem{fdef}{Definition 2.2}
\renewcommand{\thefdef}{}
\begin{fdef}
${\mathcal W}=R_2^6$ is the direct sum of six copies of $R_2$ and
$\epsilon_{m}(\overline{a})=(0,\cdots,\overline{a},\cdots,0)$ is
the element of $\mathcal W$, of which the $m$-th component is
$\overline{a}$ and others are zero, for $1\leq m\leq 6$.
\end{fdef}

Let $S_4$ be the symmetric group of $\{1,2,3,4\}$.
$$P=\{(i,j,k,l)|\{i,j,k,l\}=\{1,2,3,4\}\}$$ is the set of all the
quadruple with the  distinct components. $S_4$ has a natural
transitive action on $P$ given by
$\sigma((i,j,k,l))=(\sigma(i),\sigma(j),\sigma(k),\sigma(l))$, for
any $\sigma\in S_4$. $$H=\{(1),(13),(24),(13)(24)\}$$ is a
subgroup of $S_4$ with $[S_4:H]=6$. Then $S_4$ has a partition of
cosets with respect to $H$, denoted by
$S_4=\bigsqcup_{m=1}^6\sigma_mH$. We can obtain a partition of
$P$, $P=\bigsqcup_{m=1}^6P_m$, where $P_m=(\sigma_mH)((1,2,3,4))$.
We define the index map
$$\theta:P\rightarrow\{1,2,3,4,5,6\},$$ by
$$\theta\left((i,j,k,l)\right)=m \text{ if } (i,j,k,l)\in P_m,$$ for $1\leq
m\leq 6$. Particularly, we fix $P_1=H((1,2,3,4))$, then we have
$(1,2,3,4)\in P_1$ and $\theta((1,2,3,4))=1$.

Using the decomposition (1.8) of $\stf$, We take a $K$-basis
$\Gamma$ of $\stf$, which contains
$\{X_{ij}(r)|r\in\{r_\lambda\}_{\lambda\in\Lambda},1\leq i\neq
j\leq4\}$. Define $\psi: \Gamma \times \Gamma \to \mathcal W$ by
$$\psi(X_{ij}(r),X_{kl}(s))=\epsilon_{\theta((i,j,k,l))}(\overline{rs})\in\mathcal
W,$$ for $r,s\in\{r_\lambda\}_{\lambda\in\Lambda}$ and distinct
$i, j, k, l$ and $\psi=0$, otherwise. Then we obtain the
$K$-bilinear map $\psi:\stf\times\stf\to\mathcal W$ by linearity.

Recall that a Lie algebra over $K$ is defined to be an algebra
satisfying $[x, x]=0$ and $[[x,y], z] + [[y, z], x] + [[z, x],
y]=0$. Note the anti-commutativity does not imply $[x, x]=0$ for
arbitrary $K$.

We now have
\newtheorem{slemma}{Lemma 2.3}
\renewcommand{\theslemma}{}
\begin{slemma} The bilinear map $\psi$ is a $2$-cocycle.
\end{slemma}
{\bf Proof: }First,
$\overline{ab}=\overline{a}\overline{b}=\overline{b}
\overline{a}=\overline{ba}=-\overline{ba}$ for $a,b\in R$ and if
$(i,j,k,l)=\sigma((1,2,3,4))$($\sigma\in S_4$), then
$(k,l,i,j)=(\sigma\circ(13)(24))((1,2,3,4))$. So $(i,j,k,l)$ and
$(k,l,i,j)$ are in the same $P_m$, i.e.
$\theta((i,j,k,l))=\theta((k,l,i,j))$. Thus
$$\psi(X_{ij}(a),X_{kl}(b))=\epsilon_{\theta((i,j,k,l))}(\overline{ab})
=-\epsilon_{\theta((k,l,i,j))}(\overline{ba})=-\psi(X_{kl}(b),X_{ij}(a))$$
So $\psi$ is skew-symmetric. Since $\psi(\gamma,\gamma)=0$ for all
$\gamma\in\Gamma$, we have $\psi(x,x)=0$ for every $x\in\stf$ .

Next we prove
$$\psi([x,y],z)+\psi([y,z],x)+\psi([z,x],y)=0,$$ for any
$x,y,z\in\stf$. Denote the left side by $J(x,y,z)$, it suffices to
check $J(x,y,z)=0$ on the basis $\Gamma$. According to Lemma 1.7
and Lemma 1.11, the Steinberg Lie algebra $\stf$ has the
decomposition :
\begin{align} \stf=& t(R,R)\oplus T_{12}(1,R)\oplus T_{13}(1,R)\oplus
T_{14}(1,R)\notag \\
&\oplus_{1\leq i \neq j\leq n}X_{ij}(R), \tag{2.4}
\end{align}
where $t(R,R)$ is the $K$-linear span of the elements $t(a,b)$

We will show the following possibilities:

{\bf Case 1:} Clearly, the number of elements of ${x,y,z}$
belonging to the subalgebra $\frak T$ such that $\psi([x,y],z)\neq
0$ is at most one. Thus we can suppose that
$x=X_{12}(a),y=X_{34}(b)$ and $z\in \frak T$, where $a,b\in R$. We
omit the other subcases since they are very similar (although not
identical). By (2.2), we can assume that either $z=t(c,d)$, where
$c,d\in R$, or $z=T_{1j}(1,c)$, where $2\leq j\leq 4$ and $c\in
R$. When $z=t(a,b)$, then according to the Jacobi identity, we
have
\begin{eqnarray}
J(x,y,z))&=&\psi([t(c,d),X_{12}(a)],X_{34}(b))\nonumber\\
&=&\psi(X_{12}((cd-dc)a),X_{34}(b))\nonumber\\
&=&\epsilon_1(\overline{(cd-dc)ba)}=0;\nonumber
\end{eqnarray}
when $z=T_{12}(c)$,
\begin{eqnarray}
J(x,y,z))&=&\psi([T_{12}(1,c),X_{12}(a)],X_{3,4}(b))\nonumber\\
&=&\psi(X_{12}(ca+ac),X_{34}(b))\nonumber\\
&=&\epsilon_1(\overline{(ca+ac)b})=0;\nonumber
\end{eqnarray}
when $z=T_{13}(c)$,
\begin{eqnarray}
J(x,y,z))&=&\psi([X_{3,4}(b),T_{13}(c)],X_{12}(a))+\psi([T_{13}(1,c),X_{12}(a)],X_{3,4}(b))\nonumber\\
&=&\psi(X_{3,4}(cb),X_{12}(a))+\psi(X_{12}(ca),X_{3,4}(b))\nonumber\\
&=&\epsilon_1(-\overline{cba}+\overline{cab})=\epsilon_1(\overline{c(ab-ba)})=0;\nonumber
\end{eqnarray}
when $z=T_{14}(c)$,
\begin{eqnarray}
J(x,y,z)&=&\psi([X_{34}(b),T_{14}(c)],X_{12}(a))+\psi([T_{14}(1,c),X_{12}(a)],X_{3,4}(b))\nonumber\\
&=&\psi(-X_{34}(bc),X_{12}(a))+\psi(X_{12}(ca),X_{3,4}(b))\nonumber\\
&=&\epsilon_1(\overline{bca}+\overline{cab})=\epsilon_1(\overline{abc-abc})=0.\nonumber
\end{eqnarray}
 {\bf Case 2:} If there is none of
$\{x,y,z\}$ belonging to $\frak T$, the nonzero terms of
$J(x,y,z)$ must be $\psi([X_{ik}(a),X_{kj}(b)],X_{kl}(c))$ or
$\psi([X_{il}(a),X_{lj}(b)],X_{kl}(c))$, for distinct $i,j,k,l$
and $a,b,c\in R$.

One is: $x=X_{ik}(a)$,$y=X_{kj}(b)$ and $z=X_{kl}(c)$
\begin{eqnarray}
J(x,y,z))&=&\psi(X_{ij}(ab),X_{kl}(c))+\psi(-X_{il}(ac),X_{kj}(b))\nonumber\\
&=&\epsilon_{\theta((i,j,k,l))}(\overline{abc})-\epsilon_{\theta((i,l,k,j))}(\overline{acb})\nonumber\\
&=&\epsilon_{\theta((i,j,k,l))}(\overline{a(bc-cb)})=0.\nonumber
\end{eqnarray}

The other is: $x=X_{il}(a)$,$y=X_{lj}(b)$ and $z=X_{kl}(c)$
\begin{eqnarray}
J(x,y,z))&=&\psi(X_{ij}(ab),X_{kl}(c))+\psi(X_{il}(a),X_{kj}(cb))\nonumber\\
&=&\epsilon_{\theta((i,j,k,l))}(\overline{abc})+\epsilon_{\theta((i,l,k,j))}(\overline{acb})\nonumber\\
&=&\epsilon_{\theta((i,j,k,l))}(\overline{a(bc+cb)})=0\nonumber
\end{eqnarray}
as $(i,j,k,l)$ and $(i,l,k,j)$ are in the same partition of $P$,
i.e. $\theta((i,j,k,l))=\theta((k,l,i,j))$. This is because that
if $(ijkl)=\sigma((1234))$($\sigma\in S_4$), then
$(ilkj)=(\sigma\circ(24))((1234))$. The proof is completed.
$\Box$

We therefore obtain a central extension of Lie algebra $\stf$:
\begin{equation}
0\rightarrow{\mathcal
W}\rightarrow\stfh\overset{\pi}\rightarrow\stf\rightarrow
0,\tag{2.5}
\end{equation}
i.e. \begin{equation} \stfh={\mathcal W}\oplus\stf,\tag{2.6}
\end{equation}with Lie bracket
$$[(c,x),(c',y)]=(\psi(x,y),[x,y])$$
for all $x,y\in\stf$ and $c,c'\in{\mathcal W}$, where
$\pi:{\mathcal W}\oplus\stf\rightarrow\stf$ is the second
coordinate projection map. Then, $(\stfh,\pi)$ is a covering (or a
central extension) of $\stf$. We will show that $(\stfh,\pi)$ is
the universal covering of $\stf$. To do this, we define a Lie
algebra $\stfs$ to be the Lie algebra generated by the symbols
$X_{ij}^{\sharp}(a)$,  $a\in R$ and the $K$-linear space
${\mathcal W}$, satisfying the following relations:
\begin{align} &a\mapsto X_{ij}{^\sharp}(a) \text{ is a $K$-linear mapping,}\tag{2.7}\\
&[X_{ij}^{\sharp}(a), X_{jk}^{\sharp}(b)] = X_{ik}^{\sharp}(ab), \text{ for distinct } i, j, k, \tag{2.8}\\
&[X_{ij}^{\sharp}(a),{\mathcal W}]=0, \text{ for distinct } i, j, \tag{2.9}\\
&[X_{ij}^{\sharp}(a),X_{ij}^{\sharp}(b)]=0, \text{ for distinct } i, j, \tag{2.10}\\
&[X_{ij}^{\sharp}(a),X_{ik}^{\sharp}(b)]=0, \text{ for distinct } i, j, k, \tag{2.11}\\
&[X_{ij}^{\sharp}(a),X_{kj}^{\sharp}(b)]=0, \text{ for distinct } i, j, k, \tag{2.12}\\
&[X_{ij}^{\sharp}(a),
X_{kl}^{\sharp}(b)]=\epsilon_{\theta((i,j,k,l))}(\overline{ab}),
\text{ for distinct } j, k, i, l, \tag{2.13}
\end{align}
where $a,b\in R$,$1\leq i,j,k, l\leq 4$. As $1\in R$, $\stfs$ is
perfect. Clearly, there is a unique Lie algebra homomorphism
$\rho:\stfs\rightarrow\stfh$ such that
$\rho(X^\sharp_{ij}(a))=X_{ij}(a)$ and $\rho|_{\mathcal W}=id$.

{\bf Remark 2.14: }Comparing with the generating relations of
$\st$ (1.3)-(1.5), we separate the case $[X_{ij}^{\sharp}(a),
X_{kl}^{\sharp}(b)] ( j\neq k, i\neq l)$ into four subcases
(2.10)-(2.13).

We claim that $\rho$ is actually an isomorphism.

\newtheorem{tlemma}{Lemma 2.15}
\renewcommand{\thetlemma}{}
\begin{tlemma} $\rho:\stfs\rightarrow\stfh$ is a Lie algebra isomorphism.
\end{tlemma}
{\bf Proof: }Let
$T_{ij}^{\sharp}(a,b)=[X_{ij}^{\sharp}(a),X_{ji}^{\sharp}(b)]$.
Then one can easily check that for $a,b\in R$ and distinct
$i,j,k$, one has
\begin{align}
&T_{ij}^{\sharp}(a,b)=-T_{ji}^{\sharp}(b,a)                       \tag{2.16}\\
&T_{ij}^{\sharp}(ab,c)=T_{ik}^{\sharp}(a,bc)+T_{kj}^{\sharp}(b,ca)\tag{2.17}
\end{align}
Indeed, the proof of (2.17) is the same as the proof in [KL] and
[G1,Proposition 2.17]. Put
$t^{\sharp}(a,b)=T_{1j}^{\sharp}(a,b)-T_{1j}^{\sharp}(1,ab)$ for
$a,b\in R, 2\leq j\leq 4$. Then $t^\sharp(a,b)$ does not depend on
the choice of $j$. Also, one can easily check (as in [AF, Lemma
1.1]) that
$$\stfs={\frak T}^{\sharp}\oplus_{1\leq i \neq j\leq 4}
X_{ij}^{\sharp}(R)$$ where
$${\frak T}^{\sharp}= \left(\sum_{i, j, k, l \text{ are distinct}}[ X_{ij}^{\sharp}(R), X_{kl}^{\sharp}(R)]\right)
\oplus\left(\sum_{1\leq i < j\leq 4}[ X_{ij}^{\sharp}(R),
X_{ji}^{\sharp}(R)]\right).$$ It then follows from (2.16) and
(2.17) above that
\begin{equation}{\frak T}^{\sharp}={\mathcal W}\oplus\left(t^{\sharp}(R,R)\oplus T_{12}^{\sharp}(1,R)\oplus T_{13}^{\sharp}(1,R)\oplus
T_{14}^{\sharp}(1,R)\right)\tag{2.18}
\end{equation} where
$t^{\sharp}(R,R)$ is the linear span of the elements
$t^{\sharp}(a,b)$. So by Lemma 1.11, it suffices to show that the
restriction of $\rho$ to $t^{\sharp}(R,R)$ is injective.

Now the similar argument as  given in [AG,Lemma 6.18] shows that
 there exists a linear map from $t(R,R)$ to $t^{\sharp}(R,R)$ so
that $t(a,b)\mapsto t^{\sharp}(a,b)$ for $a,b\in R$. This map is
the inverse of the restriction of $\rho$ to $t^{\sharp}(R,R)$.
$\Box$

The following theorem is the main result of this section:
\newtheorem{thmf}{Theorem 2.19}
\renewcommand{\thethmf}{}
\begin{thmf}$(\stfh,\pi)$ is the universal covering of $\stf$ and
hence
$$H_2(\stf)\cong\mathcal W.$$
\end{thmf}
{\bf Proof: }The idea for proving this theorem is motivated by
[G1] and [AG]. Particularly, we imitate the method of proving the
universal covering of ${\frak stu}_4({\mathcal A},-,\gamma)$ in
the Section 6 of [AG].

 Suppose that
\begin{equation}0\rightarrow{\mathcal
V}\rightarrow\stft\overset{\tau}{\rightarrow}\stf\rightarrow\notag
0\end{equation} is a central extension of $\stf$. We must show
that there exists a Lie algebra homomorphism
$\eta:\stfh\rightarrow\stft$ so that $\tau\circ\eta=\pi$. Thus, by
Lemma 2.15, it suffices to show that there exists a Lie algebra
homomorphism $\xi:\stfs\rightarrow\stft$ so that
$\tau\circ\xi=\pi\circ\rho$.

Using the $K$-basis $\{r_\lambda\}_{\lambda\in\Lambda}$ of $R$, we
choose a preimage $\widetilde{X}_{ij}(a)$ of $X_{ij}(a)$ under
$\tau$, $1\leq i\neq j\leq 4,a\in
\{r_\lambda\}_{\lambda\in\Lambda}$, so that the elements
$\widetilde{X}_{ij}(a)$ satisfy the relations (2.7)-(2.13). For
distinct $i ,j, k$, let
$$[\widetilde{X}_{ik}(a),\widetilde{X}_{kj}(b)]=\widetilde{X}_{ij}(ab)+{\mu}_{ij}^k(a,b)$$
where ${\mu}_{ij}^k(a,b)\in{\mathcal V}$. Take distinct $
i,j,k,l$, then
$$\left [\widetilde{X}_{ik}(a),[\widetilde{X}_{kl}(c),\widetilde{X}_{lj}(b)]\right
]=[\widetilde{X}_{ik}(a),\widetilde{X}_{kj}(cb)].$$ But the left
side is, by Jacobi identity,
$$\left [[\widetilde{X}_{ik}(a),\widetilde{X}_{kl}(c)],\widetilde{X}_{lj}(b)]\right
]+\left
[\widetilde{X}_{kl}(c),[\widetilde{X}_{ik}(a),\widetilde{X}_{lj}(b)]\right
]=[\widetilde{X}_{il}(ac),\widetilde{X}_{lj}(b)].$$ as
$[\widetilde{X}_{ik}(a),\widetilde{X}_{lj}(b)]\in \mathcal V$.
Thus
\begin{equation}
[\widetilde{X}_{ik}(a),\widetilde{X}_{kj}(cb)]=[\widetilde{X}_{il}(ac),\widetilde{X}_{lj}(b)]\notag
\end{equation}In particular, ${\mu}_{ij}^k(a, cb)={\mu}_{ij}^l(ac, b)$ and
$[\widetilde{X}_{ik}(a),\widetilde{X}_{kj}(b)]=[\widetilde{X}_{il}(a),\widetilde{X}_{lj}(b)]$.
It follows that ${\mu}_{ij}^k(a,b)={\mu}_{ij}^l(a,b)={\mu}_{ij}(a,
b)$ which show ${\mu}_{ij}^k(a,b)$ is independent of the choice of
$k$ and ${\mu}_{ij}(c, b)={\mu}_{ij}(1, cb)$, we have
\begin{equation}
[\widetilde{X}_{ik}(a),\widetilde{X}_{kj}(b)]=\widetilde{X}_{ij}(ab)+{\mu}_{ij}(a,b)\notag
\end{equation}
Taking $a=1$, we have
\begin{equation}
[\widetilde{X}_{ik}(1),\widetilde{X}_{kj}(b)]=\widetilde{X}_{ij}(b)+{\mu}_{ij}(1,b)\notag
\end{equation} Now, we replace $\widetilde{X}_{ij}(b)$ by
$\widetilde{X}_{ij}(b)+{\mu}_{ij}(1,b)$. Then the elements
$\widetilde{X}_{ij}(b)$ still satisfy the relations (2.7).
Moreover we have
\begin{equation}
[\widetilde{X}_{ik}(a),\widetilde{X}_{kj}(b)]=\widetilde{X}_{ij}(ab)\tag{2.20}
\end{equation}
for $a,b\in R$ and distinct $i,j,k$. So, the elements
$\widetilde{X}_{ij}(a)$ satisfy (2.8).

Next for $k\neq i , k\neq j$, we have
\begin{align}
[\widetilde{X}_{ij}(a),\widetilde{X}_{ij}(b)]&=\left[\widetilde{X}_{ij}(a),[\widetilde{X}_{ik}(b),\widetilde{X}_{kj}(1)]\right]\notag\\
&=\left[[\widetilde{X}_{ij}(a),\widetilde{X}_{ik}(b)],\widetilde{X}_{kj}(1)\right]+
\left[\widetilde{X}_{ik}(b),[\widetilde{X}_{ij}(a),\widetilde{X}_{kj}(1)]\right]\notag\\
&=0+0=0\tag{2.21}
\end{align}
as both of $[\widetilde{X}_{ij}(a),\widetilde{X}_{ik}(b)]$ and
$[\widetilde{X}_{ij}(a),\widetilde{X}_{kj}(1)]$ are in $\mathcal
V$. Thus, the relation (2.10) has been shown.

For (2.11), taking $l\notin \{i,j,k\}$
\begin{align}
[\widetilde{X}_{ij}(a),\widetilde{X}_{ik}(b)]&=\left[\widetilde{X}_{ij}(a),[\widetilde{X}_{il}(b),\widetilde{X}_{lk}(1)]\right]\notag\\
&=\left[[\widetilde{X}_{ij}(a),\widetilde{X}_{il}(b)],\widetilde{X}_{kj}(1)\right]+
\left[\widetilde{X}_{il}(b),[\widetilde{X}_{ij}(a),\widetilde{X}_{lk}(1)]\right]\notag\\
&=0+0=0\tag{2.22}
\end{align}
with
$[\widetilde{X}_{ij}(a),\widetilde{X}_{il}(b)],[\widetilde{X}_{ij}(a),\widetilde{X}_{lk}(1)]\in
\mathcal V$. Similarly, we have
\begin{equation}
[\widetilde{X}_{ij}(a),\widetilde{X}_{kj}(b)]=0 \tag{2.23}
\end{equation} for distinct $i,j,k$, which is respect to the relation (2.12).

To verify (2.13)  one needs a few more steps. First, set
$\widetilde{T}_{ij}(a,b)=[\widetilde{X}_{ij}(a),\widetilde{X}_{ji}(b)]$,
The following brackets are easily checked by the Jacobi identity.
\begin{align}[\widetilde{T}_{ij}(a,b),\widetilde{X}_{ik}(c)]=\widetilde{X}_{ik}(abc), & \ \
[\widetilde{T}_{ij}(a,b),\widetilde{X}_{kj}(c)]
=\widetilde{X}_{kj}(cba)\notag\\
\text{ and }[\widetilde{T}_{ij}(a,b),\widetilde{X}_{kl}(c)]&=0
\tag{2.24}
\end{align}
then we have
\begin{align}
[\widetilde{T}_{ij}(a,b),\widetilde{X}_{ij}(c)]&=\left[\widetilde{T}_{ij}(a,b),[\widetilde{X}_{ik}(c),\widetilde{X}_{kj}(1)]\right]\notag\\
&=\left[[\widetilde{T}_{ij}(a,b),\widetilde{X}_{ik}(c)],\widetilde{X}_{kj}(1)\right]+
\left[\widetilde{X}_{ik}(c),[\widetilde{T}_{ij}(a,b),\widetilde{X}_{kj}(1)]\right]\notag\\
&=\widetilde{X}_{ij}(abc)+\widetilde{X}_{ij}(cba)=\widetilde{X}_{ij}(abc+cba)\tag{2.25}
\end{align}
for $a,b,c\in R$ and distinct $i,j,k,l$.

Next for distinct $i,j,k,l$, let
$$[\widetilde{X}_{ij}(a),\widetilde{X}_{kl}(b)]=\nu^{ij}_{kl}(a,b)$$
where $\nu^{ij}_{kl}(a,b)\in \mathcal V$. By (2.24) and (2.25),
\begin{align}
2\nu^{ij}_{kl}(a,b)&=[\widetilde{X}_{ij}(2a),\widetilde{X}_{kl}(b)]
=\left[[\widetilde{T}_{ij}(1,1),\widetilde{X}_{ij}(a)],\widetilde{X}_{kl}(b)]\right]\notag\\
&=\left[\widetilde{T}_{ij}(1,1),[\widetilde{X}_{ij}(a),\widetilde{X}_{kl}(b)]\right]+
\left[[\widetilde{T}_{ij}(1,1),\widetilde{X}_{kl}(b)],\widetilde{X}_{ij}(a)]\right]\notag\\
&=0\notag
\end{align}
which yields
\begin{equation}
\nu^{ij}_{kl}(a,b)=-\nu^{ij}_{kl}(a,b).\tag{2.26}
\end{equation}
Using Jacobi identity, we have
\begin{eqnarray}
\nu^{il}_{kj}(bc,a)&=&
[\widetilde{X}_{il}(bc),\widetilde{X}_{kj}(a)]=\left
[[\widetilde{X}_{ik}(b),\widetilde{X}_{kl}(c)],\widetilde{X}_{kj}(a)\right]\nonumber\\
&=&\left[[\widetilde{X}_{ik}(b),\widetilde{X}_{kj}(a)],\widetilde{X}_{kl}(c)\right
]+\left
[\widetilde{X}_{ik}(b),[\widetilde{X}_{kl}(c),\widetilde{X}_{kj}(a)]\right]\nonumber\\
&=&[\widetilde{X}_{ij}(ba),\widetilde{X}_{kl}(c)]=\nu^{ij}_{kl}(ba,c).\nonumber
\end{eqnarray}
It then follows that
\begin{equation}\nu^{il}_{kj}(b,a)=\nu^{ij}_{kl}(a,b)=\nu^{ij}_{kl}(ba,1)\tag{2.27}
\end{equation}where $a,b\in R$ and $i,j,k,l$ are distinct.
So $\nu^{ij}_{kl}(R,R)=\nu^{ij}_{kl}(R,1)$.

Moreover, by (2.25), we get
\begin{align}
\nu^{ij}_{kl}(abc+cba,d)&=[\widetilde{X}_{ij}(abc+cba),\widetilde{X}_{kl}(d)]=\left
[[\widetilde{T}_{ij}(a,b),\widetilde{X}_{ij}(c)],\widetilde{X}_{kl}(1)\right
]\nonumber \\
&=\left[[\widetilde{T}_{ij}(a,b),\widetilde{X}_{kl}(d)],\widetilde{X}_{ij}(c)\right
]+\left
[\widetilde{T}_{ij}(a,b),[\widetilde{X}_{ij}(c),\widetilde{X}_{kl}(d)]\right]\nonumber \\
&=[0,\widetilde{X}_{ij}(c)]+[\widetilde{T}_{ij}(a,b),\nu^{ij}_{kl}(c,d)]=0\tag{2.28}
\end{align}
Taking $c=d=1$ gives $\nu^{ij}_{kl}(ab+ba,1)=0$, i.e.
\begin{equation}\nu^{ij}_{kl}(a,b)=\nu^{ij}_{kl}(ba,1)=-\nu^{ij}_{kl}(ab,1)=\nu^{ij}_{kl}(ab,1)=\nu^{ij}_{kl}(b,a)\tag{2.29}
\end{equation}
Letting $c=1$ in (2.28) gives us
\begin{align}\nu^{ij}_{kl}(d(ab-ba),1)&=\nu^{ij}_{kl}(ab,d)-\nu^{ij}_{kl}(ba,d)=\nu^{ij}_{kl}(ab,d)+\nu^{ij}_{kl}(ba,d)\notag\\
&=\nu^{ij}_{kl}(ab+ba,d)=0\tag{2.30}
\end{align}for $a,b,c,d\in R$ and distinct $i,j,k,l$.

 By Lemma 2.1, ${\mathcal I}_2=span\{2a,c(ab-ba)|a,b,c\in R)\}$,
(2.26) and (2.30) show us
\begin{equation}\nu^{ij}_{kl}({\mathcal I}_2,1)=0.\tag{2.31}\end{equation}

Finally, recall that $P=\bigsqcup_{m=1}^6P_m$ and $\theta$ is the
index map. For any  $a\in R$, $S_4$ can also act on the set
$\{\nu^{ij}_{kl}(a,1)|(i,j,k,l)\in P\}$ by
$\sigma(\nu^{ij}_{kl}(a,1))=\nu^{i'j'}_{k'l'}(a,1)$ if
$\sigma((i,j,k,l))=(i',j',k',l')$. For each index of the partition
of $P$, we define a linear morphism
$$\xi_{(i,j,k,l)}:R\rightarrow{\mathcal
V} \ \ \text{ by } \ \ a\mapsto\nu^{ij}_{kl}(a,1).$$
 By (2.31), it induces a morphism
$R_2 \rightarrow{\mathcal V}$, which is still denoted by
$\xi_{(i,j,k,l)}$. Moreover, by (2.26),(2.27) and (2.29), we get
$$\nu^{12}_{34}(a,1)=\nu^{14}_{32}(a,1)=\nu^{34}_{12}(a,1).$$
It shows that the subgroup $H=\{(1),(13),(24),(13)(24)\}$ of $S_4$
fixes $\nu^{12}_{34}(a,1)$. Thus, we can define
$\xi_m=\xi_{(i,j,k,l)}$ if $m=\theta((i,j,k,l))$ (i.e.
$(i,j,k,l)\in P_m$), which is independent on the choice of the
elements in $P_m$. Let $\xi=\sum_{m=1}^6\xi_m\circ p_m:{\mathcal
W}\rightarrow{\mathcal V}$, where $p_m$ is the $m$-th coordinate
projection map of ${\mathcal W}=R_2^6$. Then, $\xi$ can be
extended to a Lie algebra homomorphism
$\xi:\stfs\rightarrow\stft$, by
$\xi(X^\sharp_{ij}(a))=\widetilde{X}_{ij}(a)$ for $a\in R$ and
$1\leq i\neq j\leq 4$. But then
$\tau\circ\xi(X^\sharp_{ij}(a))=\tau\circ(\widetilde{X}_{ij}(a))=X_{ij}(a)$
and $\pi\circ\rho(X^\sharp_{ij}(a))=\pi(X_{ij}(a))=X_{ij}(a)$, and
thus $\tau\circ\xi=\pi\circ\rho$ are required. We have proved this
theorem. $\Box$

\noindent {\bf Remark 2.32 }If $2$ is an invertible element of
$K$, then $R=2R$. Thus ${\mathcal I}_2=R$ and ${\mathcal
W}=R_2^6=0$. In this case, $\stf$ is centrally closed.

If the characteristic of $K$ is $2$, we display the following two
examples which are two extreme cases.

\noindent{\bf Example 2.33}  Let $R$ be an associative commutative
$K$-algebra where char $K=2$,  then we have ${\mathcal I}_2 =0$
and $R_2 = R$. Therefore $H_2(\stf)=R^6$.

\noindent{\bf Example 2.34} Let $K$ be a field of characteristic
two. $R=W_k$ is the Weyl algebra which is a unital associative
algebra over $K$ generated by $x_1, \dots, x_k, y_1, \dots, y_k$
subject to the relations $x_ix_j=x_jx_i, \ y_iy_j=y_jy_i, \
x_iy_j-y_jx_i=\delta_{ij}$. Then ${\mathcal I}_2=R$, $H_2(\stf)=0$
and $\stf$ is centrally closed.

\vspace{4mm}

\noindent{\bf \S 3 Coverings of $\stt$}

In this section we compute $H_2(\stt)$. Recall that ${\mathcal
I}_3=3R+R[R,R]$ and $R_3=R/{\mathcal I}_3$ is an associative
commutative $K$-algebra.

\newtheorem{sdef}{Definition 3.1}
\renewcommand{\thesdef}{}
\begin{sdef}
${\mathcal U}=R_3^6$ and denote the six copies of $R_3$ by
$R_3^{(1)}$,$R_3^{(2)}$,$R_3^{(3)}$ and $R_3^{(-1)}$,
$R_3^{(-2)}$,$R_3^{(-3)}$. $\overline{a}^{(i)}$ and
$\overline{a}^{(-i)}$ denote the corresponding element
$(0,\cdots,\overline{a},\cdots,0)\in {\mathcal U}$, where
$\overline{a}\in R_3$.
\end{sdef}

For convenience, we use the symbol:
$$sign(m,n)=\cases
1,&\text { if } m<n\\-1,&\text { if } m>n\endcases, \text{ for }
1\leq m\neq n\leq 3.$$
 Taking the $K$-base $\Gamma$ of $\stt$
containing $\{X_{ij}(r)\}$
$({r\in\{r_\lambda\}_{\lambda\in\Lambda},1\leq i\neq j\leq 3})$.
We define $\psi:\Gamma\times\Gamma\rightarrow\mathcal U$ by
$$\psi(X_{ij}(r),X_{ik}(s))=sign(j,k)(\overline{rs})^{(i)}$$
$$\psi(X_{ij}(r),X_{kj}(s))=sign(i,k)(\overline{rs})^{(-j)}$$ for $r,s\in
\{r_\lambda\}_{\lambda\in\Lambda}$ and distinct $i,j,k$, and
$$\psi(x,y)=0, \text { otherwise. }$$
Then $\psi$ can be extended to a bilinear map
$\stt\times\stt\rightarrow\mathcal U$.

We have
\newtheorem{plemma}{Lemma 3.2}
\renewcommand{\theplemma}{}
\begin{plemma} The bilinear map $\psi$ is a $2$-cocycle.
\end{plemma}
{\bf Proof: }By the definition, $\psi$ is skew-symmetric. Since
$\psi(\gamma,\gamma)=0$ for all $\gamma$ in the basis $\Gamma$ of
$\stt$, $\psi(x,x)=0$ for any $x\in\stt$.

Similarly to the proof of Lemma 2.3, we show  $J(x,y,z)=0$, for
$x,y,z\in\stt$. It suffices to verify this by choosing $x,y,z$
from $\Gamma$. By the decomposition,
\begin{equation}
\stt=t(R,R)\oplus T_{12}(1,R)\oplus T_{13}(1,R) \oplus_{1\leq i
\neq j\leq 3}X_{ij}(R)\tag{3.3}\end{equation}

We consider the following cases.

{\bf Case 1:} There is at most one element of $\{x,y,z\}$ in the
subalgebra $\frak T$. Let $z\in \frak T$.

We first verify two subcases $x=X_{12}(a),y=X_{13}(b)$ and
$x=X_{21}(a),y=X_{31}(b)$ for $a,b\in R$.  By (3.3), we assume
that either $z=t(c,d)$, where $c,d\in R$, or $z=T_{1j}(1,c)$,where
$2\leq j\leq 3$ and $c\in R$.

If $x=X_{12}(a),y=X_{13}(b)$, we have, according to the Jacobi
identity, when $z=t(c,d)$,
\begin{eqnarray}
J(x,y,z))&=&\psi([X_{13}(b),t(c,d)],X_{12}(a))+\psi([t(c,d),X_{12}(a)],X_{13}(b))\nonumber\\
&=&-\psi(X_{13}((cd-dc)b), X_{12}(a))+\psi(X_{12}((cd-dc)a),X_{13}(b))\nonumber\\
&=&(\overline{a(cd-dc)b+(cd-dc)ab})^{(1)}=0;\nonumber
\end{eqnarray}
when $z=T_{12}(1,c)$,
\begin{eqnarray}
J(x,y,z))&=&\psi([X_{13}(b),T_{12}(1,c)],X_{12}(a))+\psi([T_{12}(1,c),X_{12}(a)],X_{13}(b))\nonumber\\
&=&\psi(-X_{13}(cb), X_{12}(a))+\psi(X_{12}(ca+ac),X_{13}(b))\nonumber\\
&=&(\overline{acb+(ca+ac)b})^1=(\overline{3abc})^{(1)}=0;\nonumber
\end{eqnarray}
when $z=T_{13}(1,c)$,
\begin{eqnarray}
J(x,y,z))&=&\psi([X_{13}(b),T_{13}(1,c)],X_{12}(a))+\psi([T_{13}(1,c),X_{12}(a)],X_{13}(b))\nonumber\\
&=&\psi(-X_{13}(cb+bc), X_{12}(a))+\psi(X_{12}(ca),X_{13}(b))\nonumber\\
&=&(\overline{a(cb+bc)+cab})^1=(\overline{3abc})^{(1)}=0.\nonumber
\end{eqnarray}
On the other hand, if $x=X_{21}(a),y=X_{31}(b)$, when $z=t(c,d)$,
\begin{eqnarray}
J(x,y,z))&=&\psi([X_{31}(b),t(c,d)],X_{21}(a))+\psi([t(c,d),X_{21}(a)],X_{31}(b))\nonumber\\
&=&\psi(X_{31}(b(cd-dc))X_{21}(a))-\psi(X_{21}(a(cd-dc)),X_{31}(b))\nonumber\\
&=&-(\overline{(ab(cd-dc)+a(cd-dc)b)})^{(-1)}=0;\nonumber
\end{eqnarray}
when $z=T_{12}(1,c)$,
\begin{eqnarray}
J(x,y,z))&=&\psi([X_{31}(b),T_{12}(1,c)],X_{21}(a))+\psi([T_{12}(1,c),X_{21}(a)],X_{31}(b))\nonumber\\
&=&\psi(X_{31}(bc), X_{21}(a))-\psi(X_{21}(ca+ac),X_{31}(b))\nonumber\\
&=&(-\overline{abc+(ca+ac)b})_1=(-\overline{3abc})^{(-1)}=0;\nonumber
\end{eqnarray}
when $z=T_{13}(1,c)$,
\begin{eqnarray}
J(x,y,z))&=&\psi([X_{31}(b),T_{13}(1,c)],X_{21}(a))+\psi([T_{13}(1,c),X_{21}(a)],X_{31}(b))\nonumber\\
&=&\psi(X_{31}(cb+bc), X_{21}(a))-\psi(X_{21}(ac),X_{31}(b))\nonumber\\
&=&-(\overline{a(cb+bc)+acb})_1=(-\overline{3abc})^{(-1)}=0.\nonumber
\end{eqnarray}

As for the other subcases, they are similar to the above subcases
except the following subcase which does not appear above. When
$x=X_{23}(a), y=X_{21}(b)$ and $z=t(c, d)$,
\begin{eqnarray}
J(x,y,z))&=&\psi([X_{21}(b),t(c,d)],X_{23}(a))+\psi([t(c,d),X_{23}(a)],X_{21}(b))\nonumber\\
&=&\psi(X_{21}(b(cd-dc)),X_{23}(a))+\psi(0,X_{21}(b))\nonumber\\
&=&(\overline{b(cd-dc)a})^{(2)}=0.\nonumber
\end{eqnarray}

 {\bf Case 2:} If
there is none of $\{x,y,z\}$ belonging to $\frak T$, the
nontrivial terms of $J(x,y,z)$ must be
$\psi([X_{ik}(a),X_{kj}(b)],X_{ik}(c))$ or
$\psi([X_{ik}(a),X_{kj}(b)],X_{kj}(c))$, for $a,b,c\in R$ and
distinct $i,j,k$.

One is: $x=X_{ik}(a)$,$y=X_{kj}(b)$ and $z=X_{ik}(c)$
\begin{eqnarray}
J(x,y,z))&=&\psi([X_{ik}(a),X_{kj}(b)],X_{ik}(c))+\psi(X_{kj}(b),X_{ik}(c)],X_{ik}(a))\nonumber\\
&=&\psi(X_{ij}(ab),X_{ik}(c))+\psi(-X_{ij}(cb),X_{ik}(a))\nonumber\\
&=&(\overline{abc-cba})^{(i)}=0\nonumber
\end{eqnarray}

The other is: $x=X_{ik}(a)$,$y=X_{kj}(b)$ and $z=X_{kj}(c)$
\begin{eqnarray}
J(x,y,z))&=&\psi([X_{ik}(a),X_{kj}(b)],X_{kj}(c))+\psi(X_{kj}(c),X_{ik}(a)],X_{kj}(b))\nonumber\\
&=&\psi(X_{ij}(ab),X_{kj}(c))+\psi(-X_{ij}(ac),X_{kj}(b))\nonumber\\
&=&(\overline{abc-acb})^{(-j)}=0\nonumber
\end{eqnarray}

The proof is completed.    $\Box$

It is similar to the $\stf$ case, we obtain a central extension of
$\stt$.
\begin{equation}
0\rightarrow{\mathcal
U}\rightarrow\stth\overset{\pi}\rightarrow\stt\rightarrow
0,\tag{3.4}
\end{equation}
i.e. \begin{equation} \stth={\mathcal U}\oplus\stt,\tag{3.5}
\end{equation}
and define $\stts$ to be the Lie algebra generated by the symbols
$X_{ij}^{\sharp}(a)$ ,$a\in R$ and the $K$-linear space ${\mathcal
U}$, satisfying the following relations:
\begin{align} &a\mapsto X_{ij}{^\sharp}(a) \text{ is a $K$-linear mapping,}\tag{3.6}\\
&[X_{ij}^{\sharp}(a), X_{jk}^{\sharp}(b)] = X_{ik}^{\sharp}(ab), \text{ for distinct } i, j, k, \tag{3.7}\\
&[X_{ij}^{\sharp}(a),{\mathcal U}]=0, \text{ for distinct } i, j, \tag{3.8}\\
&[X_{ij}^{\sharp}(a),X_{ij}^{\sharp}(b)]=0, \text{ for distinct } i, j, \tag{3.9}\\
&[X_{ij}^{\sharp}(a),X_{ik}^{\sharp}(b)]=sign(j,k)(\overline{ab})^i, \text{ for distinct } i, j, k, \tag{3.10}\\
&[X_{ij}^{\sharp}(a),X_{kj}^{\sharp}(b)]=sign(i,k)(\overline{ab})_j,
\text{ for distinct } i, j, k, \tag{3.11}
\end{align}
where $a,b\in R$, $1\leq i,j,k\leq 3$ are distinct. $\stts$ is
perfect and there is a unique Lie algebra homomorphism
$\rho:\stts\rightarrow\stth$. We have separated the case (1.5)
into three subcases (3.9)-(3.11).

As was done in Lemma 2.15, we have
\newtheorem{llemma}{Lemma 3.12}
\renewcommand{\thellemma}{}
\begin{llemma} $\rho:\stts\rightarrow\stth$ is a Lie algebra isomorphism.
\end{llemma}

Now we can state  the main theorem of this section.
\newtheorem{fhmf}{Theorem 3.13}
\renewcommand{\thefhmf}{}
\begin{fhmf}$(\stth,\pi)$ is the universal covering of $\stt$ and
hence
$$H_2(\stt)\cong\mathcal U.$$
\end{fhmf}
{\bf Proof: } The idea to prove this theorem is similar to the
proof of Theorem 2.19. But there are some slight differences. The
point is that since $1\leq i, j, k\leq 3$, if $i, j, k$ are
distinct, then $k$ is uniquely determined once $i, j$ are chosen.

Suppose that
\begin{equation}0\rightarrow{\mathcal
V}\rightarrow
\sttt\overset{\tau}{\rightarrow}\stt\rightarrow\notag
0\end{equation} is a central extension of $\stt$. We must show
that there exists a Lie algebra homomorphism
$\eta:\stth\rightarrow\sttt$ so that $\tau\circ\eta=\pi$. Thus, by
Lemma 3.12, it suffices to show that there exists a Lie algebra
homomorphism $\xi:\stts\rightarrow\sttt$ so that
$\tau\circ\xi=\pi\circ\rho$.
 Choose a preimage $\widetilde{X}_{ij}(a)$ of $X_{ij}(a)$ as in
 Section 2.

\noindent {\bf (i).} Again let
$\widetilde{T}_{ij}(a,b)=[\widetilde{X}_{ij}(a),\widetilde{X}_{ji}(b)]$,
then
$$[\widetilde{T}_{ik}(1,1),\widetilde{X}_{ij}(a)]=\widetilde{X}_{ij}(a)+{\mu}_{ij}(a)$$
where ${\mu}_{ij}(a)\in{\mathcal V}$. Replacing
$\widetilde{X}_{ij}(a)$ by $\widetilde{X}_{ij}(a)+\mu_{ij}(a)$. By
Jacobi identity, we have
$$\left[[\widetilde{T}_{ik}(1,1),\widetilde{X}_{ik}(a)],\widetilde{X}_{kj}(b)\right]+
\left[\widetilde{X}_{ik}(a),[\widetilde{T}_{ik}(1,1),\widetilde{X}_{kj}(b)]\right]
=\left[\widetilde{T}_{ik}(1,1),[\widetilde{X}_{ik}(a),\widetilde{X}_{kj}(b)]\right]$$
which yields
$$2[\widetilde{X}_{ik}(a),\widetilde{X}_{kj}(b)]-[\widetilde{X}_{ik}(a),\widetilde{X}_{kj}(b)]
=[\widetilde{T}_{ik}(1,1),\widetilde{X}_{ij}(ab)]$$ i.e.
\begin{equation}
[\widetilde{X}_{ik}(a),\widetilde{X}_{kj}(b)]=\widetilde{X}_{ij}(ab)\tag{3.14}
\end{equation}
This gives (3.7).  The proof of the relation (3.9) is exactly same
as (2.21).

\noindent {\bf (ii).} To show $\widetilde{X}_{ij}(a)$ satisfy
(3.10) and (3.11), We define
$$[\widetilde{X}_{ij}(a),\widetilde{X}_{ik}(b)]=\nu^i_{jk}(a,b) \ \text{ and }\  [\widetilde{X}_{ij}(a),\widetilde{X}_{kj}(b)]
=\nu^{ik}_j(a,b)$$ where $\nu^i_{jk}(a,b),\nu^{ik}_j(a,b)\in
\mathcal V$. Then
\begin{align}
\nu^i_{jk}(a,b)&=[\widetilde{X}_{ij}(a),\widetilde{X}_{ik}(b)]=
\left[\widetilde{X}_{ij}(a),[\widetilde{X}_{ij}(1),\widetilde{X}_{jk}(b)]\right]\notag\\
&=[\widetilde{X}_{ij}(1),\widetilde{X}_{ik}(ab)]=\nu^i_{jk}(1,ab).\tag{3.15}
\end{align}
Moreover, we have
\begin{align}
0&=\left[\widetilde{T}_{ij}(1,1),[\widetilde{X}_{ij}(a),\widetilde{X}_{ik}(1)]\right]\notag\\
&=\left[[\widetilde{T}_{ij}(1,1),\widetilde{X}_{ij}(a)],\widetilde{X}_{ik}(1)]\right]+
\left[\widetilde{X}_{ij}(a),[\widetilde{T}_{ij}(1,1),\widetilde{X}_{ik}(1)]\right]\notag\\
&=2[\widetilde{X}_{ij}(a),\widetilde{X}_{ik}(1)]+[\widetilde{X}_{ij}(a),\widetilde{X}_{ik}(1)]=3\nu^i_{jk}(1,a)
=\nu^i_{jk}(1,3a).\tag{3.16}
\end{align}
Similarly, we have \begin{equation} \nu^{ik}_j(a,b)= \nu^{ik}_j(1,
ab) \ \text{ and } \ \nu^{ik}_j(1, 3a)=0. \tag{3.17}\end{equation}

Set
$\widetilde{t}(a,b)=\widetilde{T}_{1j}(a,b)-\widetilde{T}_{i1j}(1,ba)$,
which does not depend on the choice of $j$. By the Jacobi identity
and (3.16),
\begin{align}
0&=\left[\widetilde{t}(a,b),[\widetilde{X}_{12}(1),\widetilde{X}_{13}(c)]\right]\notag\\
&=\left[[\widetilde{t}(a,b),\widetilde{X}_{12}(1)],\widetilde{X}_{13}(c)]\right]+
\left[\widetilde{X}_{12}(1),[\widetilde{t}(a,b),\widetilde{X}_{13}(c)]\right]\notag\\
&=[\widetilde{X}_{12}(ab-ba),\widetilde{X}_{13}(c)]+[\widetilde{X}_{12}(1),\widetilde{X}_{13}((ab-ba)c)]\notag\\
&=\nu^1_{23}(1,2(ab-ba)c)=-\nu^1_{23}(1,(ab-ba)c)\notag
\end{align}
i.e. $\nu^1_{23}(1,(ab-ba)c)=0$. Also,
\begin{align}
0&=\left[\widetilde{t}(a,b),[\widetilde{X}_{12}(1),\widetilde{X}_{32}(c)]\right]\notag\\
&=\left[[\widetilde{t}(a,b),\widetilde{X}_{12}(1)],\widetilde{X}_{32}(c)]\right]\notag\\
&=[\widetilde{X}_{12}(ab-ba),\widetilde{X}_{32}(c)]\notag\\
&=\nu^{13}_{2}(ab-ba, c)=\nu^{13}_{2}(1,(ab-ba)c).\notag
\end{align}

One can show the general cases,
\begin{equation}\nu^i_{jk}(1,(ab-ba)c)=0 \text {     and    }
\nu^{ik}_j(1,(ab-ba)c)=0 \tag{3.18}
\end{equation}
for $a,b,c\in R$ and distinct $1\leq i,j,k\leq 3$.

Since ${\mathcal I}_3 = span\{3a,(ab-ba)c | a,b,c\in R\}$, by
(3.16), (3.17) and (3.18), we have
\begin{equation}\nu^i_{jk}(1, {\mathcal I}_3)=0 \text { and }
\nu^{ik}_j(1, {\mathcal I}_3)=0 \tag{3.19}
\end{equation}

Finally,
\begin{equation}
 \nu^i_{kj}(1,a)=[\widetilde{X}_{ik}(1),\widetilde{X}_{ij}(a)]
 =-[\widetilde{X}_{ij}(a),\widetilde{X}_{ik}(1)]=-\nu^i_{jk}(a,1)=-\nu^i_{jk}(1,a)\notag
\end{equation}
If we denote $\nu^i_{jk}(1,a)=\nu^i(1,a)$ for $j<k$, then
$\nu^i_{jk}(1,a)=-\nu^i(1,a)$. The same assumption can be used on
$\nu^{ik}_j(1,a)$, such that $\nu^{ik}_j(1,a)=\cases \nu_j(1,a)
&\text { for } i<k,\\-\nu_j(1,a) &\text { for } i>k.\endcases$

Taking $\xi:{\mathcal U}\rightarrow\sttt$ by
$\xi((\overline{a})^{(i)})=\nu^i(1,a)$,
$\xi((\overline{a})^{(-i)})=\nu_i(1,a)$ and $\xi|_{\sttt}=id$.

Then the rest of the proof is similar to Theorem 2.19. $\Box$

\noindent {\bf Remark 3.20}  If $R$ is commutative and char $K=3$,
 then we have ${\mathcal I}_3=3R+R[R,R]=0$, $R_3=R$. Thus,
$H_2(\stt)=R^6\neq 0$.  [KL] indicated that $H_2(\stt)=0$ when
$\frac{1}{2} \in K$. This claim is not true. Again, as in Example
2.34, if $R=W_k$ the Weyl algebra, then ${\mathcal I}_3 = R$,
$H_2(\stt)=0$ and $\stt$ is centrally closed.

\vspace{4mm}

\noindent{\bf \S 4 Concluding remarks}


Combining Theorem 1.12, Theorem 2.19 and Theorem 3.13, we
completely determined $H_2(\st)$ for $n\geq 3$.

\newtheorem{ghmf}{Theorem 4.1}
\renewcommand{\theghmf}{}
\begin{ghmf} let $K$ be a unital commutative ring  and $R$ be a
unital associative $K$-algebra. Assume that $R$ has a $K$-basis
containing the identity element. Then
\begin{equation}
H_2(\st)=\cases 0 &\text{ for } n\geq 5\\
R_2^6  &\text{ for } n=4\\
R_3^6 &\text{ for } n=3\endcases\notag
\end{equation}

\end{ghmf}

It then follows from [KL] that

\newtheorem{ihmf}{Theorem 4.2}
\renewcommand{\theihmf}{}
\begin{ihmf} let $K$ be a unital commutative ring  and $R$ be a
unital associative $K$-algebra. Assume that $R$ has a $K$-basis
containing the identity element. Then
\begin{equation}
H_2(sl_n(R))=\cases HC_1(R) &\text{ for } n\geq 5\\
R_2^6\oplus HC_1(R) &\text{ for } n=4\\
R_3^6\oplus HC_1(R) &\text{ for } n=3\endcases\notag
\end{equation}
where $HC_1(R)$ is the first cyclic homology group of the
associative $K$-algebra $R$(See [KL] or[L]).

\end{ihmf}

Department of Mathematics and Statistics

York University

Toronto, Ontario

Canada  M3J 1P3

ygao@@yorku.ca

and

Department of Mathematics

University of Science and Technology of China

Hefei, Anhui

P. R. China  230026

skshang@@mail.ustc.edu.cn

\end{document}